\documentstyle[]{article}
\author {Anna Erschler (Dyubina) \\ e-mail: annadi@math.tau.ac.il,
erschler@pdmi.ras.ru}
\title{On drift and entropy growth for random walks
on groups}
\date{preliminary version}

\begin{document}
\maketitle

\newtheorem{rem}{Remark}
\newtheorem{definition}{Defintion}
\newtheorem{lemma}{Lemma}
\newtheorem{corollary}{Corollary}
\newtheorem{proposition}{Proposition}
\newtheorem{theorem}{Theorem}
\newcommand{\erw}{{\mathop{\rm E}}}
\newcommand{\Z}{\rm{\hbox to 9pt{Z\hss Z}}}
\newcommand{\R}{{\rm R}}
\newcommand{\tgot}{{\sf t}}
\newcommand{\zz}{\Z/2\Z}

\section{Introduction}
We consider symmetric random walks on groups induced
by the measure $\mu$. In this paper we assume that the support of
$\mu$ is finite and generates the
group. We consider two functions
$$
H(n)=-\sum_{g\in G}(\mu^{*n})\ln(\mu^{*n}(g))
$$
and
$$
L(n)=\erw_{\mu^{*n}} l(g).
$$
Here $l$ denotes the word metric, corresponding to the fixed set of generators
and
$\mu^{*n}$ is the $n$-th convolution of $\mu$.
The function $H(n)$ is called the {\it entropy} and
$L(n)$ is called the {\it drift} of the random walk.
$H(n)$ measures how far the measure $\mu^{*n}$ removed from being uniformely
distributed.
$L(n)$ shows how fast (in average) the random walk is moving away
from the origin.

It is known that $H(n)$ is asymptotically linear iff
$L(n)$ is asymptotically linear (\cite{Va})
 and iff the Poisson boundary of
the random walk is nontrivial (\cite{VK}).
In particular, it is so for any non-amenable group.
On the other hand, for many examples (e.g. for
any Abelian group) $L(n)$ is asymptotically $\sqrt{n}$.

Until recently the existence problem for groups with
intermediate growth rate  of $L(n)$ was open.
First examples  were
found by the author in $\cite{Ja}, \cite{Janet}$. In these examples
$L(n)\asymp n^{1-\frac{1}{2^k}}$ (for any positive integer $k$) and
$L(n)\asymp n/\ln(n)$.

In this paper we find new possibilities for the rate of $L(n)$.
In particular, we show that  $L(n)$ can be asymptotically equal to
$$
 \frac{n}{\ln(\ln(...\ln(n)...))}.
$$
We also estimate the growth of the entropy for random walks on these groups.

The structure of the paper is the following.
In section 2 we state that certain functions
are concave. We use this auxiliary lemma in the next section.

In section 3 we consider a two-dimensional simple random walk.
We find asymptotics of some class of functionals depending on
local times of this two-dimensional random walk.

In section 4 we construct examples of groups and
we apply results of the previous section to
find asymptotics of the drift in these groups.

In section 5 we give some general estimates for the entropy $H(n)$.
We apply these estimates to the examples considered in section 4.
These examples shows that there are infinitely many possibilities
for asymptotics of the entropy. They also show that the growth of the entropy
(as well as of the drift)
can be very close to linear and yet sublinear.

\section{An auxiliary lemma}
\begin{lemma}
For
$0<\alpha<1$
let
$$
T_{k,\alpha}=\underbrace{\exp(\exp...\exp((4k)^{1/\alpha})..).}_{k}
$$

\begin{enumerate}
\item Let
$$
\tilde{L}_{k,\alpha}(x)
=\frac{x}{\underbrace{(\ln(\ln(...\ln(x))...)}_{k })^{\alpha}}.
$$
Then $\tilde{L}_{k,\alpha}(x)$ is concave on the segment $[T_{k,\alpha},\infty)$.

\item
There exists a continuous increasing
function $L_{k,\alpha}:[0,\infty) \to [0,\infty)$
and $X>0$ such that
$L_{k,\alpha}(x)$ is concave, $L_{k,\alpha}(0)=0$
 and for $x>X$ $L_{k,\alpha}(x)=\tilde{L}_{k,\alpha}(x)$

\item
 For any  $n\ge 1$ the function
$$
\frac{1}{\underbrace{\ln(\ln(...\ln(n/x))...)^\alpha }_{k}}.
$$
is concave on the segment $(0,\frac{n}{T_{k,\alpha}}]$.
\end{enumerate}
\end{lemma}
{\bf Proof.}
We prove this lemma in the Appendix.

\section{ Some functionals of two-dimensional random walk }

We say that the random walk is  simple if $\mu$ is equidistributed.
In this section we consider a simple random walk on $Z^2$.
Let $b_z^{(n)}$ be the number of times 
the random walk has 
visited the element $z$ ($z \in Z^2$) up to the moment $n$.
Let $R^{(n)}$ be the number of different elements visited until
the moment $n$. ($R^{(n)}$ is called {\it range} of the random walk).
First we formulate a simple property of the range.

\begin{lemma}
There exist $q_1, q_2>0$ such that for any $n>1$
$$
\Pr [R^{(n)} \ge q_1 \frac{n}{\ln(n)}]  \ge q_2.
$$
\end{lemma}
{\bf Proof.}
Let $\tilde{R}^{(n)}=R^{(n)} \frac{\ln(n)}{n}$.
Note that
$$
\erw [\tilde{R}^{(n)}] \asymp 1,
$$
 since (see for example \cite{Ja})
$$
\erw [R^{(n)}] \asymp \frac{n}{\ln(n)}.
$$
Hence it suffices to show that there exists $C$ such that
$$
\sigma^2 [\tilde{R}^{(n)}] \le C.
$$
In the proof of the theorem 1 \S4 \cite{Sp} it is shown that for any
random walk on a lattice
$$
\sigma^2[R^{(n})]  \le 2 \erw[R^{(n)}] \erw[R^{(n-[n/2])}+R^{(n/2)}-R^{(n)}]+
\erw[R^{(n)}] \le
6 (\erw[R^{(n)}])^2 +\erw[R^{(n)}].
$$
Hence
$$
\sigma^2 [\tilde{R}^{(n)}] \le 6 +\frac{1}{\erw[R^{(n)}]}.
$$
This completes the proof of the lemma.

\begin{lemma}
Let $T>0$ and $f$ be a concave strictly increasing function such that
\begin{enumerate}
\item For any $C>1$ $Cf(x) \ge f(Cx)$,
\item $f(0)=0$ and there exists $x$ such that $f(x)>1$,
\item For any $n\ge 1$
$xf(n/x)$ is concave on the segment $(0,n/T]$.
\end{enumerate}
Consider a simple random walk on $\Z^2$. Then there exists $K>0$ such that
$$
\erw\left(\sum_{z\in\Z^2}f(b_z^{(n)}) \right) \le
K f(\ln(n))\frac{n}{\ln(n)}.
$$
\end{lemma}

{\bf Proof.}
Note that
$$
\sum f\left(\sum b_z^{(n)} \right) \le f(n/R^{(n)})R^{(n)},
$$
since f is concave. Hence
$$
\erw\left[\sum f\left(\sum b_z^{(n)} \right)\right] \le
\erw[f(n/R^{(n)})R^{(n)}]=
$$
$$
 \erw[f(n/R^{(n)})R^{(n)} |0, R^{(n)}< n/T]+
\erw[f(n/R^{(n)})R^{(n)} | R^{(n)}\ge n/T].
$$
Let $M$ be the maximum of $f$ on the segment $[0,T]$,
that is $M=f(T)$, since $f$ is
an increasing function. Then the second
term is not greater than
$$
 \erw[M R^{(n)} | R^{(n)}\ge n/T] \le M\erw[R^{(n)}]
\asymp \frac{n}{\ln(n)}.
$$
But for $n$ large enough we have $f(\ln(n))>1$ and hence
$$
\frac{n}{\ln(n)}\le f(\ln(n))\frac{n}{\ln(n)}.
$$
Now we want to estimate the first term. Since $xf(n/x)$ is concave on
$(0,n/T]$ we have
$$
 \erw[f(n/R^{(n)})R^{(n)} |0, R^{(n)}< n/T]\le
$$
$$
\erw[R^{(n)} | 0<R^{(n)}< n/T] f\left( \frac{n}
{ \erw[R^{(n)} |0, R^{(n)}< n/T]} \right) \le
$$
$$
E[R^{(n)}] f(\frac{n}{\erw[R^{(n)}]})\asymp f(\ln(n))\frac{n}{\ln(n)}.
$$

 Then there exists $K>0$ such that
$$
\erw\left(\sum_{z\in\Z^2}f(b_z^{(n)}) \right) \le
K f(\ln(n))\frac{n}{\ln(n)}.
$$

This completes the proof of the lemma.

The following lemma gives an estimate from the other side.

\begin{lemma}
Let $f$ be a strictly increasing function on $[0,\infty)$ such that
$f(0)=0$ and for any  $C>1$ $f(Cx) \le Cf(x)$.
Then for $n$ large enough and for some positive $\epsilon$ we have
$$
\erw\left[\sum_{z\in\Z^2}f(b_z^{(n)}) \right] \ge
\epsilon f(\ln(n))\frac{n}{\ln(n)}.
$$

\end{lemma}
{\bf Proof.}
Note that for any  $\varepsilon_1 >0$ there exists $K>0$ such that
for $n$ large enough
$$
\Pr\left[b_0^{(n)} \ge K\ln(n) \right] \ge 1-\varepsilon_1
$$
(This follows from theorem 1 \cite{Erd}).
Let $n \ge 4$ and let $m=\left[ n/2 \right]$. Since $m>1$ lemma 2 implies
that
$$
\Pr [R^{(m)} \ge q_1 \frac{m}{\ln(m)}]  \ge q_2.
$$
Let $x_1^{(n)},..., x_s^{(n)}$ be different points visited by the random
walk up to the moment $n$, enumerated in the order of visiting.

Let $\beta_i^n = b_{x_i^{(n)}}^{(n)}$.

Take $\varepsilon_1$ such that $\varepsilon_1 \le 1/2 q_2$

Note that for any $0\le i \le q_1 \frac{m}{\ln(m)}$
$$
\Pr\left[\beta_i^{(n)} \ge K\ln(m) \right] \ge
\Pr[R^{(m)} \ge i]
\Pr\left[\beta_i^{(n)} \ge K\ln(m) \right]\ge 
$$

$$
q_2
\Pr\left[\beta_i^{(n)} \ge K\ln(m) |
R^{(m)} \ge i \right]  \ge
$$

$$
q_2
\Pr\left[\beta_0^{(m)} \ge K\ln(m) \right]
\Pr \left[   R^{(m)} \ge i \right]  \ge
$$

$$
q_2^2(1-\varepsilon_1) \ge q_2^2/2
$$

Hence for any  $n$ large enough

$$
\erw\left[\sum_{z\in\Z^2}f(b_z^{(n)}) \right] \ge
\frac{q_2^2}{2} f(K\ln(\left[ n/2 \right] )) q_1 \frac{n/2}{\ln(n/2)} \asymp
f(\ln(n))\frac{n}{\ln(n)}.
$$


\begin{corollary}
\begin{enumerate}
\item Let $L_{k,\alpha}(x)$ be the function defined in Lemma 1
($0<\alpha<1$).
Then
$$
 \erw\left[\sum_{z\in\Z^2}L_{k,\alpha}(b_z^{(n)})\right] \asymp
L_{k+1,\alpha}(n).
$$

\item
 Let $f(x)=x^\alpha (0<1<\alpha)$ then
$$
 \erw\left[\sum_{z\in\Z^2}f(b_z^{(n)}) \right]
\asymp n/\ln(n)^{(1-\alpha)}
$$
\end{enumerate}
\end{corollary}
{\bf Proof.}
This corollary follows from lemma 2 and lemma 3 since for $n>N$
$$
\frac{n}{\ln(n)} L_{k,\alpha}(\ln(n)) \asymp
\frac{n}{\ln(n)} \frac{\ln(n)}{(\ln(\ln...\ln(n)...))^{\alpha}} \asymp
L_{k+1,\alpha}(n)
$$
and
$$
\frac{n}{\ln(n)} \ln(n)^{1-\alpha}= n/\ln(n)^\alpha.
$$

\section{Main result}

First we recall the definition of the wreath product.

\begin{definition}
The wreath product of  $C$ and $D$ is a semidirect product
$C$ and $\sum_C D$, where $C$ acts on
$\sum_C D$ by shifts:
if $c\in C$, $f:C \to D, f\in \sum_C D$,  then
$f^c(x)=f(xc^{-1}), x\in C$.
Let $C\wr D$ denote the wreath product.
\end{definition}

\begin{lemma}
Let $a_1, a_2,..., a_k$ generate A.
Consider the symmetric
simple random walk on $A$ corresponding to this set of generators.
Let
$$
L(n)=L_n^A.
$$
Then for some simple random walk on $B=\Z^2\wr A$
$$
L_n^B\asymp \erw \sum_{z\in \Z^2}L(b_z^{(n)})
$$
\end{lemma}
{\bf Proof.}
The proof of this lemma is similar to that of lemma 3 in \cite{Janet}.
 For any $a\in A$
$\tilde{a}^e$ denotes the function  from $\Z^2$ to $A$ such that
$\tilde{a}^e(0)=a$ and $\tilde{a}^e(x)=e$ for any $x\ne 0$. Let $a^e=
(e,\tilde{a}^e)$.
Let $e_1', e_2'$ be the standard generators of $\Z^2$ and
 $e_1=(e_1',e)$,  $e_2=(e_2',e)$.

Consider the following set of generators of $B$:
$$(a_j^e)^p e_s (a_n^e)^q,$$
 $p,q=0,1$ or $-1$,
 $s=1$ or  $2$ and
$1\le j,n \le k$.

Consider the simple random walk on $B$, corresponding to this set
of generators.

Suppose that this random walk hits  $(a,f)$ at the moment $n$. Let
$l_A(f(z))=c_z$.
Then
 $$ \frac{1}{2} \sum c_z \le l(a,f)\le 2(\sum c_z+R)$$
To see this note that multiplying by one of the generators changes the
value of $f$ in at most 2 points and that two multiplyings suffice
if  we want to change the value
of $f$ in $a$ standing at the point $a$.

Hence

$$
\frac{1}{2}\erw \left[\sum c_z \right]
=\frac{1}{2}\sum\erw[c_z] \le \erw [l(a,f)] \le 2(\sum\erw[c_z]+\erw[R])
$$

The projection of this random walk on $\Z$ is simple and symmetric.
Note that in any $i\in \Z^2$ there is a simple random walk on $A$.
The measure that defines this random walk is supported on
 $a_i$($1\le i \le n$), $e$, and inverses of these elements.
It is uniformely distributed on
$a_i$ and their inverses.
Let $\hat{L}_n^A$ be the drift that corresponds to this new measure.
It is clear that
$$
\hat{L}_n^A \asymp L_n^A.
$$


Note that

$$\erw[c_z | b_z=b, z\ne 0, a\ne z]=\hat{L}_{2b}^A$$
$$\erw[c_z | b_z=b, z=0, a\ne z]=\hat{L}_{2b-1}^A$$
$$\erw[c_z | b_z=b, z\ne 0, a=z]=\hat{L}_{2b-1}^A$$
$$\erw[c_z | b_z=b, z=0, a=z]=\hat{L}_{2b-2}^A$$

Hence
$$
\erw\left[\min(\hat{L}_{2b}^A,
 \hat{L}_{2b-1}^A, \hat{L}_{2b-2})\right] \le
\erw[c_z|b_z=b]\le
 \erw\left[\max(\hat{L}_{2b}^A,
\hat{L}_{2b-1}^A, \hat{L}_{2b-2})\right]
$$
There exist $C_2,C_3>0$ such that
$$
C_2 L(n)\le \hat{L}_{2n-2}^A,
 \hat{L}_{2n-1}^A, \hat{L}_{2n}^A \le C_3 L(n)
$$

Hence
$$
C_2 \erw[\sum_{i\in Z^2}L(b_z)] \le \erw[c_z] \le C_3 \erw[\sum_{i\in Z^2}L(b_z)]
$$
This completes the proof of the lemma.

\begin{theorem}
\begin{enumerate}
\item Let $F$ be a finite group.
Consider the following groups  that are defined recurrently
$$G_1=\Z^2\wr F; G_{i+1}=\Z^2\wr G_i.$$
Then for some simple random walk on $G_i$ and for any $n$ large enough
we have
$$
L_n^{G_i} \asymp \frac{n}{\underbrace{\ln(\ln...\ln{n})...)}_k}
$$
\item
Consider the following groups that are defined recurrently
$$
F_1=\Z; F_{i+1}=\Z\wr F_i
$$
and let
$$
H_{1,i}=\Z^2\wr F_i;   H_{j+1,i}=\Z^2 \wr H_{j,i}
$$
Then for some simple random walk on $H_{j,i}$ and for any $n$ large enough
we have
$$
L_{H_{j,i}}(n) \asymp \frac{n}{ \underbrace{\sqrt[2^i]{\ln(\ln(...\ln}}_j
(n)...)) }.
$$

\end{enumerate}
\end{theorem}

{\bf Proof.}
\begin{enumerate}
\item
 We prove the theorem by induction on $i$.
Base $i=1$. In this case $G_i=\Z^2\wr F$ and
$L(n)$ is asymptorically equal to $\frac{n}{\ln(n)}$ (\cite{Ja})
Induction step follows from the previous lemma and corollary 1.

\item
We prove the statement by induction on $j$.
 For $H_{1,i}=F_i$ the asymptotics of the drift is found in \cite{Janet}.
Is is proven there that for some random walk on $F_i$
$$
L_{F_i}(n) \asymp n^{1-\frac{1}{2^k}}.
$$
The induction step follows from the previous lemma and colollary 1.

\end{enumerate}

\section{Estimates of the entropy}

In this section we give estimates for the entropy of a random walk.
It is known (see \cite{der}) that for a wide class of measures on nilpotent
groups $H(n) \asymp  \ln(n)$. As it was mentioned before 
 $H(n)$ is asymptotically linear
for any nonamenable group. In this section
we study intermediate examples.

Let
$v(n)$ be the growth function of the group
(that is $v(n)=\# \{g\in G: l(g)\le n\}$).
Let
$$
l=\lim_{n\to \infty} L(n)/n,
$$
$$
h=\lim_{n\to \infty} H(n)/n,
$$
$$
v=\lim_{n\to \infty} \ln(v(n))/n.
$$
(See for example \cite{Ve} for the proof that these limits exist.)
It is known (see \cite{Ve}) that
$$
h\le v l.
$$

The following lemma generalizes this fact.
\begin{lemma}
For any $\varepsilon >0$ there exists $C>0$ such that
$$
H(n)\le (v+\varepsilon)L(n)+\ln(n)+C.
$$
\end{lemma}
{\bf Proof.} Let $a_i^{(n)}=\Pr_{\mu^{*n}} l(g)$. Then by definition
$$
L(n)=\sum_{i=0}^n  i a_i^{(n)}.
$$

Note that
$$
H(n)\le \sum_{i=1}^n a_i^{(n)} \ln (v(i)/a_i^{(n)})=
$$
$$
\sum_{i=1}^n a_i^{(n)} \ln (v(i))+\sum_{i=1}^n a_i^{(n)}(- ln (a_i^{(n)}))\le
\sum_{i=1}^n a_i^{(n)} \ln (v(i))+\ln(n).
$$
For any $\varepsilon>0$ there exists $K>0$ such that
$v(i)\le K(v+\varepsilon)^i$. Hence
$$
H(n)\le \sum_{i=1}^n a_i^{(n)} (i(v+\varepsilon)+\ln(K))+\ln(n) =
$$
$$
(v+\varepsilon)\sum_{i=1}^n a_i^{(n)} i+\ln(K))+\ln(n)=
(v+\varepsilon)L(n)+\ln(K))+\ln(n).
$$

Another lemma estimates the entropy from the other side:
\begin{lemma}
\begin{enumerate}
\item
There exists $C>0$ such that
$$
H(n)\ge C\erw_{\mu^{*n}} l^2(g)/n -\ln(n) \ge CL^2(n)/n -\ln(n)
$$
\item
There exists $K>0$ such that
$$
L(n) \le K \sqrt{n (\ln(v(n))+ \ln(n))}.
$$
\end{enumerate}
\end{lemma}

{\bf Proof.}
\begin{enumerate}
\item
Let $p_n(x)$ be the probability to hit $x$ after $n$ steps.
In \cite{Va} it is shown that there exist $K_1, K_2>0$ such
that for any $x$ and $n$
$$
p_n(x)\le K_1 n^{3/4} \exp(-K_2 l(x)^2/n) \le
 K_1 n \exp(-K_2 l(x)^2/n).
$$
Note that then
$$
-\sum_{g\in G: l(g)=i} \mu^{*n}(g) \ln(\mu^{*n}(g)) \ge
(-\ln( K_1) -\ln(n) + K_2 i^2/n)a_i^{(n)}
$$
Hence
$$
H(n) \ge -\ln(K_1) -\ln(n)+ K_2\sum_{i=0}^n i^2/n a_i^{(n)}=
$$
$$
-\ln(K_1) -\ln(n) + K_2/n\erw_{\mu^{*n}} l^2(g) \ge
C/n\erw_{\mu^{*n}} l^2(g) -\ln(n).
$$
The last inequality follow from the fact that for some $C_2>0$
$H(n)\ge C_2$.

\item This follows from the first part of the lemma, since
$H(n) \le \ln (v(n))$.
\end{enumerate}

As a corollary from the two previous lemmas we get the following theorem.
\begin{theorem}
Let $G_i$ be the groups defined in theorem 1.
Then for some random walk on $G_i$ we have
$$
K_1 n/\underbrace{\ln(\ln...\ln}_i (n)...)^2 \le
H_{G_i}(n)\le
K_2 n/\underbrace{\ln(\ln...\ln}_i (n)...)
$$
for some positive constansts $K_1$ and $K_2$.
In particular, all $G_i$ have different asymptotics of the entropy.
\end{theorem}

\appendix
\section{Proof of the auxiliary lemma}
In the appendix we give the proof of lemma 1.

\begin{enumerate}
\item
Let $m_{k,\alpha}(x)=(\underbrace{\ln(\ln...\ln(x)...)}_{k})^{\alpha}$.
We want to prove that
$x/m_{k,\alpha}(x)$ is concave on $[T_{k,\alpha},\infty)$.
Note that
$$
(x/m_{k,\alpha})''=\frac{-xm_{k,\alpha}''
m_{k,\alpha}^2-2m_{k,\alpha}'m_{k,\alpha}^2+2xm_{k,\alpha}(m_{k,\alpha}')^2}
{m_{k,\alpha}^4}.
$$
Since $m_{k,\alpha}(x)>0$ on $(0,n/T_{k,\alpha}]$ it suffices to prove that
$$
2x(m_{k,\alpha}')^2-xm_{k,\alpha}m_{k,\alpha}''\le 2m_{k,\alpha}'m_{k,\alpha}.
$$
To prove this we will show that
$$
2x(m_{k,\alpha}')^2< 1/2m_{k,\alpha}m_{k,\alpha}'
$$
and that
$$
m_{k,\alpha}(-m_{k,\alpha}'')x\le 1.5 m_{k,\alpha}'m_{k,\alpha}.
$$
Note that
$$
m_k'(x)=\frac{\alpha}
{x\underbrace{\ln(x) \ln(\ln(x))...\ln(\ln(...\ln(x)...))}
_{k-1}m_{k,\alpha}^{1-\alpha}}.
$$
We see that $0<m_{k,\alpha}'(x)<1/x$ and we know that $m_{k,\alpha}(x)>4$
($x\in [T_{k,\alpha},\infty)$). This proves the first inequality.

Now let
$$
r_{k,\alpha}(x)=\alpha/m_k'=
x\underbrace{\ln(x) \ln(\ln(x))...\ln(\ln(...\ln(x)...))}
_{k-1}
m_{k,\alpha}^{1-\alpha}.
$$

Note that
$$
r_{k,\alpha}'(x)=
m_{k,\alpha}^{1-\alpha}m_{1,1}(x)m_{2,1}(x)...m_{k-1,1}(x)
+
$$
$$
\sum_{i=1}^{k-1}r_{k,\alpha}(x)\frac{m_{i,1}'(x)}{m_{i,1}(x)}+
r_{k,\alpha}\frac{(1-\alpha)m_{k,\alpha}'}
{m_{k,\alpha}}.
$$
Note that for $1\le i \le k-1$
$m_{i,1}'(x)\le 1/x$ and that
$m_{k,alpha}'(x)\le 1/x$.
Since
$m_{k,\alpha}(x)\ge 2k$
 and for $1\le i \le k-1$
$m_i(x)\ge 2k$
we get
$$
r_{k,\alpha}'x\le 1.5r_{k,\alpha}.
$$

This implies that
$$
(-m_{k,\alpha}'')x=
\alpha xr_{k,\alpha}'/r_{k,\alpha}^2 \le
 \alpha 1.5/r_{k,\alpha} =1.5m_{k,\alpha}'.
$$
So we have proven also the second inequality.
\item
Let $\beta_k=\tilde{L}_{k,\alpha}(T_{k,\alpha})/(2T_{k,\alpha})$.
 Consider the function
$y_k(x)=\beta_k x$. Note
that $y_(T_{k,\alpha}) < \tilde{L}(T_{k,\alpha})$ and that for $x>X$
 $y_(x) > \tilde{L}(x)$.
Take maximal $z$ such that
 $y_(z)= \tilde{L}(z)$.
Let $L_{k,\alpha}(x)=y(x)$ if $0\le x \le z$ and
$L_{k,\alpha}(x)=\tilde{L}_{k,\alpha}(x)$ if
$x\ge z$. Note that $y(x)$ is concave ($x\in [0,x]$), $\tilde{L}_k(x)$
is concave if $x>z$ and that $\tilde{L}_{k,\alpha}'(z) \le y'(z)$.
This implies that $L_{k,\alpha}$ is concave.

\item
We can assume without loss of generality that $\alpha=1$
since $x^\alpha$ is a concave function.
Let $g(x)=\underbrace{ln(ln...(ln(n/x)))}_k$.
 We want to proof that $1/g$ is concave on
the given segment.
Note that
$$
(1/g)''=\frac{-g''g^2+2g(g')^2}{g^4}.
$$
Since $g$ is positive on the given segment, we need to prove that
$$
2(g')^2v\le g''g.
$$
Let
$$
h(x)=x\underbrace{\ln(n/x) \ln(\ln(n/x))... \ln(\ln...\ln(n/x)...))}_{k-1}.
$$
Note that $g'(x)=-1/h$ and $g''(x)=h'/(h^2)$.
We want to prove that $2<gh'$. Since $g(x)>2$ on the given segment, it
suffices to note, that $h'(x)>1$. This follows from the fact that
$h(x)/x$ increases and greater than $1$.
\end{enumerate}

The author expresses her gratitude to A.M.Vershik for stating
the problem and many useful discussions.

\end{document}